\documentclass[a4paper]{amsart}

\usepackage{amssymb, amsmath, amsthm, microtype, hyperref}
\hypersetup{colorlinks, urlcolor=black, citecolor=blue}

\newcommand{\R}{\mathbb{R}}
\newcommand{\N}{\mathbb{N}}

\newcommand{\vol}{\operatorname{vol}}
\renewcommand{\epsilon}{\varepsilon}
\newcommand{\PSH}{\operatorname{PSH}}
\newcommand{\nn}{\text{nn}}
\newcommand{\nK}{\text{nK}}
\newcommand{\reg}{{\operatorname{reg}}}
\newcommand{\supp}{\operatorname{supp}}
\newcommand{\inv}{^{-1}}

\newtheorem{theorem}{Theorem}[section]
\newtheorem{prop}[theorem]{Proposition}
\theoremstyle{definition}
\newtheorem{definition}[theorem]{Definition}
\newtheorem{lemma}[theorem]{Lemma}
\newtheorem{remark}[theorem]{Remark}

\title{Definitions of the volume of a big cohomology class}
\author{Tiernan Cartwright}

\begin{document}
\begin{abstract}
    We elaborate on how two definitions of the volume of a big cohomology class are consistent. The first definition involves taking the absolutely continuous part of a closed positive current, and the second involves the non-pluripolar product. We also describe how a similar equality holds for the numerical restricted volume introduced by Collins and Tosatti.
\end{abstract}
\maketitle
\section{Introduction}
Given a holomorphic line bundle $L$ over a compact K\"ahler $n$-fold $X$, an important quantity is its volume, which measures the rate of growth of the amount of global sections of tensor powers of $L$. This is defined by
\[\vol(L) = \limsup_{k\to\infty} \frac{n!}{k^n}h^0(X, kL),\]
and it measures how big a line bundle is in the sense that $\vol(L) >0$ if and only if $L$ is big. For a survey on big line bundles and their volume (mostly in the algebraic setting), see \cite{ELMNP}. For ample line bundles, from the asympototic Riemann--Roch formula \cite[1.1.25]{Positivity} we can deduce the topological interpretation that volume is the top-degree self-intersection number, i.e.\
\[\vol(L) = \int_X c_1(L)^n.\]
This explains the name `volume'. In fact, this formula still holds if $L$ is merely nef, due to the holomorphic Morse inequalities; see for instance the lecture notes \cite{PluripotentialTheory} for an exposition.

If $L$ is an arbitrary big (or, more generally, pseudo-effective) line bundle, \cite[Theorem 1.2]{Boucksom} gives an explicit analytic formula for $\vol(L)$. Namely,  
\begin{equation}\label{eq: bigL}
    \vol(L) = \max_T\int_X T^n_{ac}
\end{equation}
where $T$ ranges among the (real) closed positive $(1, 1)$-currents in the first Chern class $c_1(L)$, and $T_{ac}$ denotes the absolutely continuous part of $T$ (see Subsection \ref{ssec: lebesgue} below).  This gives a formula for $\vol(L)$ in full generality, since we recall that if $L$ is not big, $\vol(L) = 0$. In particular, the formula shows that $\vol(L)$ is an invariant of the cohomology class $c_1(L)$. This motivates the following definition.
\begin{definition}{\cite[Definition 1.3]{Boucksom}}\label{def: Boucksom} Let $X$ be a compact K\"ahler $n$-fold. We define the volume of a cohomology class $\alpha\in H^{1, 1}(X, \R)$ by 
	\[\vol(\alpha) = \sup_T \int_X T^n_{ac}\]
	for $T$ ranging over the closed positive $(1, 1)$-currents in $\alpha$, in the case $\alpha$ is pseudo-effective (which is equivalent to the existence of such currents in $\alpha$, by definition). Otherwise, we set $\vol(\alpha) = 0$.
\end{definition}
The supremum is finite \cite[Proposition 2.6]{Boucksom}, and Equation \eqref{eq: bigL} states that $\vol(c_1(L)) = \vol(L)$ for any holomorphic line bundle $L$. 

The volume of a cohomology class of a compact K\"ahler $n$-fold has a similar geometric and topological interpretation to the volume of a holomorphic line bundle over such a manifold. Specifically, $\vol(\alpha) >0$ if and only if $\alpha$ is a big class \cite[Theorem 4.7]{Boucksom}. Moreover, if $\alpha \in H^{1, 1}(X, \R)$ is nef (not necessarily big, though this is the case of primary interest), then its volume is given by the cup product
\[\vol(\alpha) = \alpha^n,\]
which is Theorem 4.1 in \cite{Boucksom}. 

As part of creating a framework for studying the complex Monge--Amp\`ere equation in a big cohomology class, the paper \cite{BEGZ} introduces a generalization of the wedge product called the non-pluripolar product, which works for arbitrary closed positive $(1, 1)$-currents. During the study of its properties, they introduce an alternative definition of the volume of a big cohomology class, which involves taking a top-degree non-pluripolar product of a current with itself; this definition will be presented later as Definition \ref{def: BEGZ}.

The purpose of this note is to show how the more recent definition of the volume of a big cohomology class in the paper \cite{BEGZ} is consistent with the original definition due to Boucksom \cite{Boucksom} stated as Definition \ref{def: Boucksom} above. This elaborates on the comment made before Proposition 1.18 in \cite{BEGZ}, where the authors indicate that it will follow from that proposition that the two definitions are consistent. Although the results in this paper may already be known to experts, we think it is worthwhile to give an exposition of them, especially since there are technicalities involved in any notion of the product of closed positive currents. We also show in Section \ref{sec: restricted} how these techniques extend to the notion of numerical restricted volume introduced by Collins and Tosatti \cite{RestrictedVol}.

\vspace*{2mm}
\noindent\textbf{Proof outline:} The proof of the main result consists of 3 ideas.
\begin{enumerate}
    \item Arbitrary currents can be approximated within their cohomology class by currents with analytic singularities (this is Demailly's regularization theorem, see Theorem \ref{thm: approx}).
    \item For currents with analytic singularities, the non-pluripolar product coincides with the wedge product of the absolutely continuous part (see Proposition \ref{prop: npp=ac}).
    \item Finally, there is a monotonicity formula (see Theorem \ref{thm:monotone}) which states that less singular currents have greater mass.
\end{enumerate}
\vspace*{2mm}
\textbf{Organization of the paper:} In Section \ref{sec:2} we introduce the necessary preliminaries, including the non-pluripolar product. In Section \ref{sec:3} we state the definition of volume due to \cite{BEGZ}, and show that the two definitions of volume coincide, which is the main result of this paper. A similar result holds in the setting of restricted volumes, which we introduce and discuss in Section \ref{sec: restricted}.

\vspace*{2mm}
\noindent\textbf{Acknowledgments:} I am very thankful to S\l awomir Dinew for his opinion and feedback on an earlier version of this manuscript. I also thank my supervisor Zhou Zhang for his support and encouragement to write this paper. This research is supported by an Australian Government Research Training Program (RTP) Scholarship.

\section{Preliminaries}\label{sec:2}
\subsection{Non-pluripolar product}
Throughout this paper, except where otherwise indicated, $X$ denotes a compact K\"ahler $n$-fold with K\"ahler form $\omega$. Plurisubharmonic functions will be called `psh' for short.

We recall the definition of a big cohomology class. The map from $H^{1, 1}(X, \R)$ to the space of closed real $(1, 1)$-currents modulo $dd^c$-exact currents, given on representatives by taking a form to its associated current, is an isomorphism. Therefore we can view a cohomology class in $H^{1, 1}(X, \R)$ as a collection of currents.
\begin{definition}\label{def: big}
	A cohomology class $\alpha\in H^{1, 1}(X, \R)$ is called pseudo-effective if it contains a closed positive $(1, 1)$-current. It is called big if it contains a K\"ahler current, i.e.\ a current $T$ such that $T\ge \varepsilon_0 \omega$ (in the sense of currents) for some constant $\varepsilon_0 >0$.
\end{definition}
To study the Monge--Amp\`ere equation in cohomology classes which are big but not necessarily K\"ahler, we need to deal with plurisubharmonic functions that may not be bounded. The non-pluripolar product gives a framework for extending the Monge--Amp\`ere operator to such functions. This was introduced in the influential paper \cite{BEGZ}, which is the primary reference for this subsection.

Let $(X, \omega)$ be an arbitrary Hermitian manifold (i.e.\ $\omega$ is a smooth strictly positive $(1, 1)$-form), and let $u_1, \dots, u_p$ be psh functions on $X$. For all $k\in \N$, let $O_k = \bigcap_{j=1}^p \{u_j > -k\}$. Then the current 
\[1_{O_k}\bigwedge_{j=1}^p dd^c\max\{u_j, -k\}\]
is non-decreasing in $k$; an analogous result for $\omega$-psh functions is proved in \cite{qpshenergy} but we will give a proof here for the reader's convenience.

For all $l \in \N$, $j=1, \dots, p$, let $\max\{u_j, -l\} = u_{j,l}$, which is a bounded psh function approximating $u_j$. If $l \ge k$, then \[\bigcap_{j=1}^p\{u_{j,l} > -k\} = \bigcap_{j=1}^p\{u_j > -k\} = O_k.\] Thus by the locality of the Bedford--Taylor product in the plurifine topology (see \cite{plurifine}), for which $O_k$ is open,
\[1_{O_k}\bigwedge_{j=1}^p dd^c\max\{u_{j,l}, -k\} = 1_{O_k}\bigwedge_{j=1}^p dd^c u_{j,l}.\]
Moreover, $l\ge k$ implies $\max\{u_{j, l},-k\} = \max\{u_j, -k\}$ and $O_k \subset O_l$. Thus
\[1_{O_k}\bigwedge_{j=1}^p dd^c\max\{u_j, -k\} \le 1_{O_l}\bigwedge_{j=1}^p dd^c\max\{u_j, -l\},\]
as claimed.

It follows that if the integral of \[1_{O_k}\bigwedge_{j=1}^p dd^c\max\{u_j, -k\}\] against any test $(n-p,n-p)$-form can be bounded independently of $k$, this sequence of currents will have a limit as $k\to\infty$.
\begin{definition}{\cite{BEGZ}} Let $(X, \omega)$ be an arbitrary Hermitian manifold. Let $u_1, \dots, u_p$ be psh functions on $X$, and for all $k\in \N$ let $O_k = \bigcap_{j=1}^p \{u_j > -k\}$. We say that the non-pluripolar product $\langle \bigwedge_{j=1}^p dd^c u_j\rangle$ is well-defined on $X$ if for each compact subset $K$ of $X$,
\[\sup_{k\in \N} \int_{K\cap O_k} \omega^{n-p} \wedge \bigwedge_{j=1}^p dd^c\max\{u_j, -k\} < \infty.\]
If it is well-defined, we define the non-pluripolar product by
\[\left< \bigwedge_{j=1}^p dd^c u_j\right> = \lim_{k\to\infty}1_{O_k}\bigwedge_{j=1}^p dd^c\max\{u_j, -k\},\]
where the limit is in the sense of currents.
\end{definition}
This definition is non-trivial only when $X$ is non-compact, since a compact manifold has no non-constant psh functions. However, the product $\langle\bigwedge_{j=1}^p dd^c u_j\rangle$ and whether it is well-defined depends only on the currents $dd^c u_j$, not the choice of potentials $u_j$ \cite[Proposition 1.4]{BEGZ}. Thus on any complex manifold $X$, we can define the non-pluripolar product of closed positive $(1, 1)$-currents $T_1, \dots, T_p$ on $X$ via taking local potentials; this is meaningful regardless of whether $X$ is compact.
\begin{definition}
We say that the non-pluripolar product $\langle T_1 \wedge \cdots \wedge T_p\rangle$ of closed positive currents is well-defined if $X$ admits an open cover $\{U_\alpha\}_{\alpha\in A}$, such that on each $U_\alpha$ there exists psh functions $u_{j,\alpha}$, such that \[T_j = dd^cu_{j, \alpha}\]
for all $j=1, \dots, p$ and $\langle \bigwedge_{j=1}^p dd^c u_{j, \alpha}\rangle$ is well-defined on $U_\alpha$.
\end{definition} In this case, $\langle T_1 \wedge \dots\wedge T_p\rangle$ is a closed positive $(p, p)$-current \cite[Theorem 1.8]{BEGZ}. Another important property of the non-pluripolar product is that on any pluripolar set it coincides with the zero current.

One reason we focus on the compact K\"ahler setting in this paper is given by the following result.
\begin{prop}[\protect{\cite[Proposition 1.6]{BEGZ}}]
   Let $X$ be compact K\"ahler. If $T_1, \dots, T_p$ are arbitrary closed positive $(1, 1)$-currents on $X$, their non-pluripolar product $\langle T_1 \wedge \dots \wedge T_p\rangle$ is well-defined. 
\end{prop}
 The compactness assumption is necessary: see \cite[Example 1.3]{BEGZ} for a counterexample in the local setting. The K\"ahler assumption can be relaxed to the assumption that $X$ is of Fujiki class $\mathcal{C}$, but it is not known if this result holds for arbitrary compact complex $X$ (see \cite[Remark 1.7]{BEGZ}).

\subsection{Mildly singular currents}\label{sec: singcurrents}
Many of the difficulties of studying $(1, 1)$-currents in a big cohomology class come from their singularities. In order to get a handle on this, there are various notions of currents which have singularities which are mild in some sense. 

In this subsection (and throughout the following), $X$ is a compact K\"ahler $n$-fold. Fix a big cohomology class $\alpha \in H^{1, 1}(X, \R)$ and a representative $\theta$. A function $\varphi: X \to [-\infty, \infty)$ is called $\theta$-plurisubharmonic ($\theta$-psh for short) if it is locally the sum of a psh function and a smooth function, and $\theta + dd^c\varphi \ge 0$ in the sense of currents. Note that $\theta$-psh functions are upper semi-continuous and thus bounded from above on $X$ (by compactness). The set of $\theta$-psh functions on $X$ is denoted $\PSH(X, \theta)$.

For any closed $(1, 1)$-current $T$ in $\alpha$, there is a $(0, 0)$-current $\varphi$ such that $T = \theta + dd^c\varphi$, called the global potential of $T$. If $T$ is positive, we can identify $\varphi$ with a $\theta$-psh function. We control the singularities of $T$ by controlling those of $\varphi$.

Let $\varphi$, $\psi \in \PSH(X, \theta)$. We say that $\varphi$ is less singular than $\psi$ if there is a constant $C$ such that
\[\psi \le \varphi + C\]
on $X$. If this holds and moreover $\psi$ is also less singular than $\varphi$, we say $\varphi$ and $\psi$ have the same singularity type, which is an equivalence relation. There is a unique equivalence class of $\theta$-psh functions which are less singular than any other $\theta$-psh function on $X$; a representative is the envelope
\[V_\theta = \sup \{\varphi \in \PSH(X, \theta) :\; \varphi \le 0\}.\]
If $T = \theta + dd^c\varphi$ and $S = \theta + dd^c\psi$ are closed positive currents, we say that $T$ is less singular than $S$ if the global potential $\varphi$ is less singular than $\psi$.
\begin{definition}
	A $\theta$-psh function is said to have minimal singularities if it has the same singularity type as $V_\theta$. A closed positive $(1, 1)$-current $T = \theta + dd^c\varphi$ is said to have minimal singularities (within its cohomology class) if its global potential $\varphi$ does.
\end{definition}

In the proof of the equivalence of the definitions of volumes, a crucial step is approximation by currents with analytic singularities.
\begin{definition}
	A closed positive $(1, 1)$-current $T = \theta + dd^c\varphi$ and its global potential $\varphi$ are said to have analytic singularities if there exists $c > 0$ such that locally on~$X$,
	\[
	    \varphi = \frac{c}{2}\log\biggl(\sum_{j=1}^N |f_j|^2\biggr) + u,
	\]
	where $u$ is smooth and $f_1, \dots, f_N$ are local holomorphic functions. 
\end{definition}
 Approximation by such currents, with only an arbitrarily small loss of positivity, is possible by a theorem due to Demailly.
 \begin{theorem}[Demailly's regularization theorem \cite{DemReg}]\label{thm: approx}
     Let $T$ be a closed $(1, 1)$-current on $(X, \omega)$ such that $T \ge \gamma$, for a smooth real $(1, 1)$-form $\gamma.$ Then there exists a sequence $T_k$ of currents with analytic singularities such that:
     \begin{enumerate}
     \item $T_k$ is cohomologous to $T$,
        \item $T_k \to T$,
        \item  $T_k \ge \gamma - \epsilon_k \omega$, where $\epsilon_k > 0$ is a sequence converging to zero, and
        \item  the Lelong numbers $\nu(T_k, x)$ increase to $\nu(T, x)$ uniformly with respect to $x\in X$.
     \end{enumerate}
 \end{theorem}
\subsection{Lebesgue decomposition}\label{ssec: lebesgue}
Given a (real) positive $(p, p)$-current $T$ written locally as
\[i^{p^2}\sum_{I, J} T_{IJ}dz_{I}\wedge d\bar{z}_J,\]
the coefficients $T_{IJ}$ can be considered as Radon measures (see \cite[Proposition 1.5]{Kolo}). Lebesgue's decomposition theorem then implies that each $T_{IJ}$ decomposes as a sum of two parts, the first which is absolutely continuous with respect to the Lebesgue measure, and the second which is singular. This gives a decomposition of $T$ itself, $T = T_{ac} + T_{\text{sing}}$, where we call $T_{ac}$ the absolutely continuous part of $T$. Note that $T_{ac}$ is in general not closed, even when $T$ is. The current $T_{ac}$ can be seen as a positive form with $L^1_{\text{loc}}$ coefficients, so we can take its exterior power to get a positive Borel $(p, p)$-form $T_{ac}^p$. We emphasise that in this notation, the $ac$ part should be taken before the $p$-th power.

\begin{lemma}\label{lem:leb}
If $T$ is a closed positive $(1,1)$-current with analytic singularities along $A\subset X$, then $T_{ac} =1_{X\setminus A}T$. Moreover, $T_{ac}$ is a closed positive $(1, 1)$-current.
\end{lemma}
This is proved by showing that the Lebesgue decomposition of a closed positive current with analytic singularities coincides with its Siu decomposition. See \cite[Subsection 2.3]{Boucksom} for details.
\section{Equivalence of definitions}\label{sec:3}
Fix a big cohomology class $\alpha \in H^{1, 1}(X, \R)$. In Definition \ref{def: Boucksom} of $\vol(\alpha)$, there is a supremum. In the non-pluripolar product definition of volume, this will correspond to choosing a current with minimal singularities, due to the following result.
\begin{theorem}\label{thm:monotone}
   Let $T$ and $S$ be cohomologous closed positive currents. If $S$ is less singular than $T$, then
\begin{equation*}
   \int_X \langle T^n\rangle \le \int_X \langle S^n\rangle. 
\end{equation*}
\end{theorem}
This was proved under the hypothesis of small unbounded loci in Theorem 1.16 of \cite{BEGZ}. This hypothesis is satisfied both for currents with analytic singularities and for currents with minimal singularities, which will suffice for our purposes. However, Theorem \ref{thm:monotone} was proved in full generality by Witt Nystr\"om in \cite{MonotonicityWN}. See also \cite{MonotonicityDDNL} for an extension to mixed~products.

We can now present the definition of volume due to \cite{BEGZ}. 
\begin{definition}\cite[Definition 1.17]{BEGZ}\label{def: BEGZ}
	Let $\alpha \in H^{1, 1}(X, \R)$ be big. Let $T_{\text{min}}$ be a positive current with minimal singularities in $\alpha$. Then the volume of~$\alpha$~is
	\[\int_X \langle T_{\text{min}}^n\rangle,\]
	which we will denote by $\langle \alpha^n\rangle$ to distinguish it from Definition \ref{def: Boucksom}.
\end{definition}
It follows from Theorem \ref{thm:monotone} that this definition is independent of the choice of current $T_{\min}$ with minimal singularities, and that if $T\in \alpha$ is any current, we~have
\begin{equation}\label{eq:min}
  \int_X\langle T^n\rangle \le \langle \alpha^n\rangle.  
\end{equation}
\begin{remark}
    Although in Definition \ref{def: BEGZ} the authors choose to only call $\langle \alpha^n\rangle$ a volume when $\alpha$ is big, it is not any less general. As remarked in \cite[page 219]{BEGZ}, the volume $\langle \alpha^n\rangle$ is continuous on the big cone, so we can extend it by continuity to the pseudo-effective cone, and set it to zero when $\alpha$ is not pseudo-effective. Since $\vol(\alpha): H^{1, 1}(X, \R) \to \R$ is continuous \cite[Corollary 4.11]{Boucksom} and is zero when $\alpha$ is not pseudo-effective (by definition), the two definitions $\vol(\alpha)$ and $\langle \alpha^n\rangle$, where the latter is in this generalised sense, will coincide for arbitrary $\alpha \in H^{1, 1}(X, \R)$ if and only if they coincide whenever $\alpha$ is big. This justifies focusing on the big case.
\end{remark}
The notation in Definition \ref{def: BEGZ} is suggestive. Once we prove that $\langle \alpha^n\rangle = \vol(\alpha)$, it will follow that when $\alpha$ is moreover nef, $\langle\alpha^n\rangle = \alpha^n$, where the latter is the cup product.

An important observation is that for currents with analytic singularities, the non-pluripolar product coincides with the absolutely continuous part. We will give a proof of this result since it involves two different generalizations of the wedge product. Note that a positive $(n,n)$-current is a positive linear functional on the space of continuous functions $C^0(X, \R).$ Via the Riesz representation theorem, it is thus associated with a Radon measure. This is association is often left implicit in the pluripotential theory literature (and it will be implicit in Section \ref{sec: restricted}), but for this proof the notions will be distinguished for clarity.
\begin{prop}\label{prop: npp=ac}
	Let $T$ be a closed positive $(1, 1)$-current on $X$ with analytic singularities. Then
	\[\langle T^n\rangle = T_{ac}^n.\]
\end{prop}
\begin{proof}
	Let $\mu$ and $\nu$ be the Radon measures associated to the positive $(n, n)$-currents $\langle T^n\rangle$ and $T_{ac}^n$ respectively. We will prove the proposition by showing that $\mu(U) = \nu(U)$ for all open $U\subset X$.

Let $U$ be an arbitrary open subset of $X$. Let $A$ be the analytic set where $T$ is singular. In particular, $A$ is pluripolar, so
\[\mu(U) = \mu(U \setminus A)\]
(since $\langle T^n\rangle$ does not charge pluripolar sets) and $A$ has Lebesgue measure zero so \[\nu(U) = \nu(U\setminus A)\] (since $\nu$ is absolutely continuous with respect to the Lebesgue measure).

By construction, the $\mu$-measure of the open set $U\setminus A$ is 
\[\mu(U \setminus A) = \sup\left\{\int_X \langle T^n\rangle \wedge f\;:\: f\in C^0(X, \R), \; 0\le f\le 1,\; \supp f \subset U\setminus A\right\},\]
and $\nu(U\setminus A)$ can be given similarly (see the proof of Theorem 2.14 in \cite{Rudin}). Fix $f\in C^0(X, \R)$ such that $0\le f\le 1$ and $\supp f \subset U\setminus A$. We calculate that
\begin{align*}
    \int_X \langle T^n\rangle \wedge f &= \int_{X\setminus A} T^n \wedge f\\
    &= \int_{X\setminus A} T^n_{ac} \wedge f\\
    &= \int_X T^n_{ac} \wedge f.
\end{align*}
This is because $T|_{X\setminus A} = T_{ac}|_{X\setminus A}$ is smooth, and for smooth currents the non-pluripolar product coincides with the wedge product. Taking the supremum over all such $f$, we conclude that $\mu(U\setminus A) = \nu(U\setminus A)$ and hence $\mu(U) = \nu(U)$.
\end{proof}

The main result of this note now follows by approximating the volumes by currents with analytic singularities.
\begin{theorem}\label{thm:main}
	Let $X$ be a compact K\"ahler $n$-fold and let $\alpha \in H^{1, 1}(X, \R)$ be a big cohomology class. Then the two definitions of volume coincide, i.e.\
	\[\vol(\alpha) = \langle \alpha^n\rangle.\]
\end{theorem}
\begin{proof}
From \cite[Proposition 1.18]{BEGZ}, there is a sequence $T_k$ of K\"ahler currents in $\alpha$ with analytic singularities such that
\[\lim_{k\to\infty}\int_X \langle T^n_{k} \rangle = \int_X \langle T^n_{\text{min}}\rangle = \langle \alpha^n\rangle,\]
where $T_{\min}$ is a choice of current with minimal singularities in $\alpha$. (Explicitly, since $\alpha$ is big it contains a K\"ahler current, and by Theorem \ref{thm: approx} it contains a K\"ahler current $T_+$ with analytic singularities. Let $S_k$ be the sequence of currents converging to $T_{\min}$ from Demailly's regularization theorem, satisfying $S_k \ge -\varepsilon_k\omega$. Then $T_k = (1-\varepsilon_k)S_k + \varepsilon_k T_+$.) Combining this with Proposition~\ref{prop: npp=ac},
\begin{align*}
	\langle \alpha^n\rangle &= \lim_{k\to\infty}\int_X \langle T^n_{k} \rangle\\
	&= \lim_{k\to\infty} \int_X T^n_{k, ac}\label{eq:npp=ac}\\
	&\le \sup_{k\in \N} \int_X T^n_{k,ac}\\
	&\le \sup_T \int_X T^n_{ac}\\
    &=\vol(\alpha).
\end{align*}

For the reverse inequality, now let $T_k$ be a sequence of positive currents  with analytic singularities such that \[\lim_{k\to\infty} \int_X T^n_{k, ac} = \sup_T \int_X T^n_{ac} = \vol(\alpha).\]
See \cite[Lemma 4.4]{Boucksom} for the existence and a description of this sequence; it also involves Demailly's regularization theorem. For each $k \in \N$, we have
\[\int_X T^n_{k, ac} = \int_X \langle T^n_{k} \rangle \le \langle\alpha^n\rangle\]
by Proposition \ref{prop: npp=ac} and Equation \eqref{eq:min}. Taking $k\to\infty$, we get
\[\vol(\alpha) \le \langle \alpha^n\rangle,\]
completing the proof.
\end{proof}
\section{Restricted volume}\label{sec: restricted}
In this section, we comment on how a similar observation can be made in the more general setting of restricted volumes. As before, $X$ is a compact K\"ahler $n$-fold. In this section, $\alpha = [\theta]$ is a pseudo-effective cohomology class.

For any closed positive (1, 1)-current $T$, let $E_+(T) = \{x \in X:\, \nu(T, x) > 0\}$ be the subset of $X$ where $T$ has positive Lelong number. Then the non-K\"ahler locus of $X$ (with respect to $\alpha$) is defined by
\[E_{\nK}(\alpha) = \bigcap_{\text{K\"ahler currents } T \in [\alpha]} E_+(T),\] and the non-nef locus is defined by
\[E_{\nn}(\alpha) = \bigcup_{\epsilon > 0} E_{\nK}(\alpha+\epsilon\omega).\]
If $\alpha$ is psuedoeffective but not big, the intersection in $E_\nK(\alpha)$ is empty, and we use the convention $E_\nK(\alpha) = X$. We always have $E_{\nn}(\alpha) \subset E_{\nK}(\alpha)$. These subsets were introduced in \cite{DivisorialZariski}.

In \cite{RestrictedVol}, Collins and Tosatti introduce the numerical restricted volume $\langle\alpha^k\rangle_{X|V}$ of $\alpha$ on the irreducible $k$-dimensional analytic subvariety $V$ of $X$, which if $V \nsubseteq E_{\nn}(\alpha)$ is given by
\[\langle\alpha^k\rangle_{X|V} = \lim_{\epsilon \to 0^+}\sup_T \int_{V_{\reg}} ((T+\epsilon\omega)|_{V_\reg})^k_{ac},\]
where the supremum is over all closed real $(1, 1)$-currents $T \in \alpha$ such that $T\ge -\epsilon\omega$ and $T$ has analytic singularities which do not contain $V$, and is equal to $0$ if $V \subseteq E_{\nn}(\alpha)$. This is a finite number \cite[Lemma 2.1]{RestrictedVol}. Notice the similarity of this definition to Definition \ref{def: Boucksom}. This extends a notion of restricted volume that had previously been defined independently by Hisamoto \cite{Bergman} and Matsumura \cite{RestrictedZariski} for the case where $V \nsubseteq E_{\nK}(\alpha)$.

The reason we discuss this notion in this note is because it can also can be recast using the non-pluripolar product of a current with minimal singularities. First we need a monotonicity formula that holds for subvarieties. 
\begin{prop}\label{prop: r monotone}
    Let $V$ be an irreducible $k$-dimensional analytic subvariety of $X$. Let $S$ and $T$ be cohomologous closed positive currents on $X$ such that $S$ is less singular than $T$. Then 
    \[\int_{V_\reg} \langle(T|_{V_\reg})^k\rangle \le \int_{V_\reg} \langle(S|_{V_\reg})^k\rangle\]
\end{prop}
\begin{proof}
    Let $f: \Tilde{X} \to X$ be an embedded resolution of singularities of $V$, i.e.\ such that $\Tilde{X}$ and the proper transform $\Tilde{V}$ are compact K\"ahler manifolds and $f|_{\Tilde{V}}: \Tilde{V} \to V$ is bimeromorphic. Let $U$ be the Zariski open subset of $V$ such that $f|_{\Tilde{U}}: \Tilde{U} \to U$ is a biholomorphism, where $\Tilde{U} = f\inv(U)$. Over $U$, we have
    \begin{align*}
        \int_{\Tilde{U}} f^*\langle T^k\rangle|_{\Tilde{U}} &= \int_{\Tilde{U}} f^*\langle (T|_U)^k\rangle\\
        &= \int_U\langle (T|_U)^k\rangle\\
        &= \int_{V_\reg}\langle (T|_{V_\reg})^k\rangle
    \end{align*} 
    where the last equality is because the Zariski closed set $V\setminus U$ is pluripolar, hence contributes no mass. Moreover,
       \[\int_{\Tilde{U}} f^*\langle T^k\rangle|_{\Tilde{U}}  = \int_{\Tilde{U}}\langle (f^*T|_{\Tilde{U}})^k\rangle = \int_{\Tilde{V}}\langle (f^*T|_{\Tilde{V}})^k\rangle.\]
       (The first equality is the biholomorphic invariance of the non-pluripolar product. This can been seen by taking local potentials as in the definition of the product; for a psh function $u$ defined in an open subset of $U$ and $j\in \N$, we have
       \begin{align*}
           f^*1_{\{u > -j\}}(dd^cu)^k &= f^*1_{\{u > -j\}} (dd^c\max\{u, -j\})^k\\ &= 1_{\{f^*u > -j\}}(dd^c\max\{f^*u, -j\})^k
       \end{align*}
       because $f^*dd^cv = dd^c(f^*v)$ for locally bounded psh $v$.)
       In summary, we have
    \begin{equation}\label{eq: pullback}
        \int_{V_\reg}\langle (T|_{V_\reg})^k\rangle = \int_{\Tilde{V}}\langle (f^*T|_{\Tilde{V}})^k\rangle,
    \end{equation}
    and the same equality holds replacing $T$ with $S$. It is easy to check that $f^*S$ is less singular than, and cohomologous to, $f^*T$. Since $\Tilde{V}$ is compact K\"ahler, we can use the usual monotonicity inequality (Theorem \ref{thm:monotone}) to get
    \[\int_{\Tilde{V}}\langle (f^* T|_{\Tilde{V}})^k\rangle \le \int_{\Tilde{V}}\langle (f^* S|_{\Tilde{V}})^k\rangle,\]
    Combining this with Equation \eqref{eq: pullback}, we conclude that
    \[\int_{V_\reg} \langle(T|_{V_\reg})^k\rangle \le \int_{V_\reg} \langle(S|_{V_\reg})^k\rangle.\qedhere\]
\end{proof}
We now show how the non-pluripolar product relates to the numerical restricted volume. In Remark 2.2 and Remark 5.2 of \cite{RestrictedVol}, the authors indicate that the following result can be proved by adjusting the proof of their Lemma 5.1. We provide a proof in order to show that it can be proved by adapting the ideas in the proof of Theorem \ref{thm:main}.
\begin{theorem}\label{thm: rvol}
   If $V$ is an irreducible $k$-dimensional analytic subvariety of the compact K\"ahler manifold $(X, \omega)$ such that $V\nsubseteq E_{\operatorname{nn}}(\alpha)$, then
\[\langle\alpha^k\rangle_{X|V} = \lim_{\epsilon \to 0^+} \int_{V_\reg} \langle (T_{\min, \epsilon}|_{V_\reg})^k\rangle,\]
where $T_{\min, \epsilon}$ is a positive current with minimal singularities in the class $[\theta + \epsilon\omega]$, where $\alpha = [\theta]$.
\end{theorem}
\begin{proof}
     We first recall the beginning of the proof of \cite[Lemma 5.1]{RestrictedVol}. For each $\epsilon > 0$, applying Demailly's regularization theorem to $T_{\min, \frac{\epsilon}{2}} - \frac{\epsilon}{2}\omega$, there exists a current $T_\epsilon \in \alpha$ with analytic singularities such that $T_\epsilon \ge -\epsilon\omega$. Since the Lelong numbers increase along Demailly's regularization, \[E_+(T_\epsilon) \subseteq E_+(T_{\min, \frac{\epsilon}{2}}) \subseteq E_\nK\left(\alpha + \frac{\epsilon}{2}\omega\right) \subseteq E_\nn(\alpha).\]
    For the second containment, we used that a current on $X$ with minimal singularities in $\alpha + \frac{\epsilon}{2}\omega$ has locally bounded potentials on the set $X \setminus E_\nK\left(\alpha + \frac{\epsilon}{2}\omega\right)$, which is non-empty since $\alpha + \frac{\epsilon}{2}\omega$ is big.

    We now describe the adjustments needed to prove Theorem \ref{thm: rvol}. Fix $\epsilon > 0$. Let $S_\epsilon$ be an arbitrary current in $\alpha$ with analytic singularities not containing $V$ and such that $S_\epsilon \ge -\frac{\epsilon}{2}\omega$. By Proposition \ref{prop: npp=ac},
    \[
        \int_{V_\reg}\left(S_\epsilon|_{V_\reg} + \frac{\epsilon}{2}\omega\right)_{ac}^k = \int_{V_\reg}\Bigl\langle \left(S_\epsilon|_{V_\reg} + \frac{\epsilon}{2}\omega\right)^k\Bigr\rangle \le \int_{V_\reg}\langle (S_\epsilon|_{V_\reg} + \epsilon\omega)^k\rangle.
    \] 
    Since $S_\epsilon + \frac{\epsilon}{2}\omega$ is a positive current in $[\theta + \frac{\epsilon}{2}\omega]$ and multiples of $\omega$ do not affect the singularity type, the current $S_\epsilon + \epsilon\omega$ is more singular than $T_{\min, \frac{\epsilon}{2}} + \frac{\epsilon}{2}\omega$, which is more singular than $T_\epsilon + \epsilon\omega$ due to the regularization. Thus $T_\epsilon + \epsilon\omega$ is less singular than $S_\epsilon + \epsilon\omega$ and cohomologous to it. Proposition \ref{prop: r monotone} then gives 
    \[\int_{V_\reg}\langle (S_\epsilon|_{V_\reg} + \epsilon\omega)^k\rangle \le \int_{V_\reg}\langle (T_\epsilon|_{V_\reg} + \epsilon\omega)^k\rangle.\] Taking the supremum over such $S_\epsilon$, and using that the limit in the definition of $\langle\alpha^k\rangle_{X|V}$ is non-increasing, we have shown 
    \[\langle\alpha^k\rangle_{X|V} \le \int_{V_\reg}\langle (T_\epsilon|_{V_\reg} + \epsilon\omega)^k\rangle.\]
    Applying Proposition \ref{prop: r monotone} twice more, we have
\begin{align*}
    \int_{V_\reg}\langle (T_\epsilon|_{V_\reg} + \epsilon\omega)^k\rangle \le& \int_{V_\reg}\langle (T_{\min, \epsilon}|_{V_\reg})^k\rangle\\
    \le& \int_{V_\reg}\langle (T_{\min, \epsilon}|_{V_\reg} +\epsilon\omega)^k\rangle\\
    \le& \int_{V_\reg}\langle (T_{2\epsilon}|_{V_\reg} + 2\epsilon\omega)^k\rangle
\end{align*}
Note that $T_{2\epsilon} \ge -2\epsilon\omega$, and $E_+(T_{2\epsilon}) \subseteq E_{\nn}(\alpha)$ implies $V \nsubseteq E_+(T_{2\epsilon})$. Thus 
\[\int_{V_\reg}\langle (T_{2\epsilon}|_{V_\reg} + 2\epsilon\omega)^k\rangle \le \sup_{R_\epsilon}\int_{V_\reg}\langle (R_{\epsilon}|_{V_\reg} + 2\epsilon\omega)^k\rangle = \sup_{R_\epsilon}\int_{V_\reg} (R_{\epsilon}|_{V_\reg} + 2\epsilon\omega)^k_{ac},\]
where the supremum is over currents in $\alpha$ satisfying $R_\epsilon \ge -2\epsilon\omega$ with analytic singularities not containing $V$. In summary, for all $\epsilon >0$, we have
\begin{equation}\label{eq:summary}
    \langle \alpha^k\rangle_{X|V} \le \int_{V_\reg}\langle (T_{\min, \epsilon}|_{V_\reg})^k\rangle \le \sup_{R_\epsilon}\int_{V_\reg} (R_{\epsilon}|_{V_\reg} + 2\epsilon\omega)^k_{ac}.
\end{equation}
By definition,
\[\lim_{\epsilon \to 0^+}\sup_{R_\epsilon}\int_{V_\reg} (R_{\epsilon}|_{V_\reg} + 2\epsilon\omega)^k_{ac} = \langle \alpha^k\rangle_{X|V},\]
so taking $\epsilon \to 0$ in Equation \eqref{eq:summary} gives the desired formula.
\end{proof}
Compared to the proof of Theorem \ref{thm:main}, this proof is similar in that it involves Proposition \ref{prop: npp=ac} and a monotonicity inequality. Note that the numerical restricted volume is already defined using currents with analytic singularities. However, Demailly's regularization theorem is still used to approximate $T_{\min, \epsilon}$. 
\bibliographystyle{alpha}
\bibliography{vol_paper.bib}
\end{document}